\documentclass[10pt]{article}
\usepackage{hyperref}
\usepackage{fancyvrb}  
\usepackage{amsmath,amssymb,color,bbm,ifthen,enumerate}
\usepackage{graphicx}
\usepackage{subfigure}
\usepackage{multirow}

\newtheorem{thm}{Theorem}[section]

\newcommand{\xv}{\mbox{\boldmath$x$}}

\newcommand{\bone}{\mbox{\boldmath$1$}}

\begin{document}

 \title{Technique for computing the PDFs and CDFs of non-negative
   infinitely divisible random variables.}

\author{Mark S. Veillette and Murad S. Taqqu \thanks{ This work was partially supported
    by the NSF grant DMS-0706786 at Boston University.} \thanks{{\em
      AMS Subject classification}. Primary 60E07,  65C50, 6008   Secondary
    6004 }
\thanks{{\em Keywords and phrases:}  Infinitely divisible
  distributions, Post-Widder formula, Stable distributions,
  Stochastic integration  }  }

\maketitle
\begin{abstract}
We present a method for computing the PDF and CDF of a
non-negative infinitely divisible random variable $X$.   Our method uses the
L\'{e}vy-Khintchine representation of the Laplace transform
$\mathbb{E} e^{-\lambda X} = e^{-\phi(\lambda)}$, where $\phi$ is the
Laplace exponent.   We apply the Post-Widder method for
Laplace transform inversion combined with a sequence convergence
accelerator to obtain accurate results.    We demonstrate this technique
on several examples including the stable distribution, mixtures
thereof, and integrals with respect to non-negative L\'{e}vy
processes.  Software written to implement this method is available from the authors and we illustrate its use at the end of the paper.
\end{abstract}

\section{Introduction}

Let $X$ be a non-negative random variable.  The distribution of $X$ is said to be infinitely divisible (ID) if for any positive integer $n$, we can find i.i.d.~random variables $X_{i,n}$, $i=1,2,\dots,n$ such that
\begin{equation*} 
X \overset{d}{=} X_{1,n} + X_{2,n} + \dots + X_{n,n}.
\end{equation*}
For background, see \cite{steutel:2004}, \cite{Applebaum:2004}.  There are many examples of such distributions, including the gamma distribution, compound Poisson distributions, inverse Gaussian distribution and right-skewed stable distributions.  These distributions are also central in the study of non-decreasing L\'{e}vy processes.   In this paper, we give a method for numerically computing the probability density function (PDF) and cumulative distribution function (CDF) of a non-negative infinitely divisible random variable.  

Our starting point is the L\'{e}vy-Khintchine (LK) formula (\cite{Bertoin:1996}, \cite{Bertoin:1999},\cite{sato:1999}), which in the case of non-negative ID random variables states that the Laplace transform of $X$,
\begin{equation}\label{e:LK1}
\psi(\lambda) \equiv \mathbb{E} e^{-\lambda X} = e^{-\phi(\lambda)}, \quad \lambda > 0, 
\end{equation}
has an exponent $\phi(\lambda)$ which is called the Laplace exponent, which can be written as
\begin{equation} \label{e:LK2}
\phi(\lambda) = a \lambda + \int_0^\infty (1 - e^{-\lambda u}) \Pi(du).
\end{equation}
  Here, $a \geq 0$ is a shift and $\Pi$ is a measure on $(0,\infty)$
  which satisfies 
\begin{equation}\label{e:intassum}
\int_0^\infty (1 \wedge x) \Pi(dx) < \infty.
\end{equation}  
Since we are interested primarily in PDFs and CDFs, we will assume
that $a = 0$ throughout, since a positive $a$ only shifts the
PDF/CDF.   Let $f_X(x)$ and $F_X(x)$ denote the PDF and CDF of $X$,
respectively (we will promptly drop the subscript $X$ on these
functions when it is clear from the context).  The Laplace transforms
of the PDF and CDF can be easily obtained from the LK formula (we will use a tilde to denote a Laplace transform):
\begin{equation*}
\psi(\lambda) = \widetilde{f}(\lambda) = \int_0^\infty e^{-\lambda x} f(x) dx, 
\end{equation*}
and a simple application of Fubini's theorem implies
\begin{equation}\label{e:derPSI}
\Psi(\lambda) \equiv \widetilde{F}(\lambda) = \int_0^\infty e^{-\lambda x} \left( \int_{0}^x f(y) dy \right) dx = \frac{1}{\lambda} \int_0^\infty e^{-\lambda y} f(y)  dy 
= \frac{\psi(\lambda) }{\lambda} .
\end{equation}

Thus, obtaining $f$ and $F$ is a matter of inverting a Laplace transform.  Generally, this task is not easy.  Typically it is done by complex integration of the Laplace transform (see \cite{abate:1992}), which can be difficult if the integrands are slowly-decaying, oscillatory functions.  This causes many numerical integration methods to converge slowly.  Here, we apply a different method of Laplace inversion known as the {\it Post-Widder (PW) method} (\cite{Grassmann:2000}, Theorem 2 or \cite{feller:1971}, section
VII.6).  It is based on the fact that under weak conditions on a
function $f$, we have  
\begin{equation}\label{e:PWformula}
f(x) = \lim_{k \rightarrow \infty} \frac{(-1)^{k-1}}{(k-1)!} \left( \frac{k}{x} \right)^{k} \tilde{f}^{(k-1)}\left(\frac{k}{x}\right), \quad x > 0,
\end{equation} 
where $\tilde{f}^{(k-1)}(k/x)$ denotes the $(k-1)^{th}$ derivative of
$\tilde{f}$ evaluated at $k/x$.  Thus, instead integrating the Laplace
transform, we are taking arbitrarily high derivatives.  

The obvious challenge in using this method is computing high
derivatives of $\tilde{f}$ in (\ref{e:PWformula}) to approximate the
limit.  Some methods have been developed to do this in general (see
for instance, \cite{jagerman:1982} or \cite{abate:1992}), however most
involve complex integration.  The method we describe here for the case
of infinitely divisible distributions is both easy to implement, and
only involves (at worst) integration of a real valued non-negative
exponentially decaying function.  The method we use here combines both
computation of a finite number of terms of the sequence in
(\ref{e:PWformula}) and an numerical extrapolation method to
approximate the limit,  \cite{frolov:1998}.  A similar method was used
to compute mean first-passage times of L\'{e}vy subordinators in
\cite{veillette1:2008}.

This paper is organized as follows.  In Section \ref{s:PWmet} we give
an overview of the approximation method used based on the (PW)
formula.  In Section \ref{s:implementation}, we outline the algorithm
for finding $f$ and $F$, as well as numerical issues that may arise in
the computation.  We test the method in cases where the PDF and CDF
are known in closed form in Section \ref{s:testit}.  In Section
\ref{s:examples} we apply our method to a collection of examples.  The
software for implementing the methods described here is freely
available form the authors and its use it described in Section \ref{s:guidetosoft}.

\section{Post-Widder method with extrapolation}\label{s:PWmet}

Given a continuous function $g(x), \ x > 0$, which is bounded as $x \rightarrow \infty$, we will denote its $k^{th}$ PW approximation as
\begin{equation}\label{e:fk}
g_k(x) = \frac{(-1)^{k-1}}{(k-1)!} \left( \frac{k}{x} \right)^{k} \tilde{g}^{(k-1)}\left(\frac{k}{x}\right), \quad x > 0.
\end{equation}
The fact that $g_k(x) \rightarrow g(x)$ as $k\rightarrow \infty$  can
be seen by approximating $g(x)$ by $\mathbb{E} g(\bar{Y}_k)$, where
$\bar{Y}_k = k^{-1}(Y_1 + Y_2 + \dots + Y_k)$, and the $Y_i$'s are
i.i.d.~gamma random variables with mean $x$ and variance $x^2$ and
then applying the law of large numbers (see Section VII.6 in \cite{feller:1971}).  

\begin{figure}[h]
\centering
\includegraphics[scale=.6]{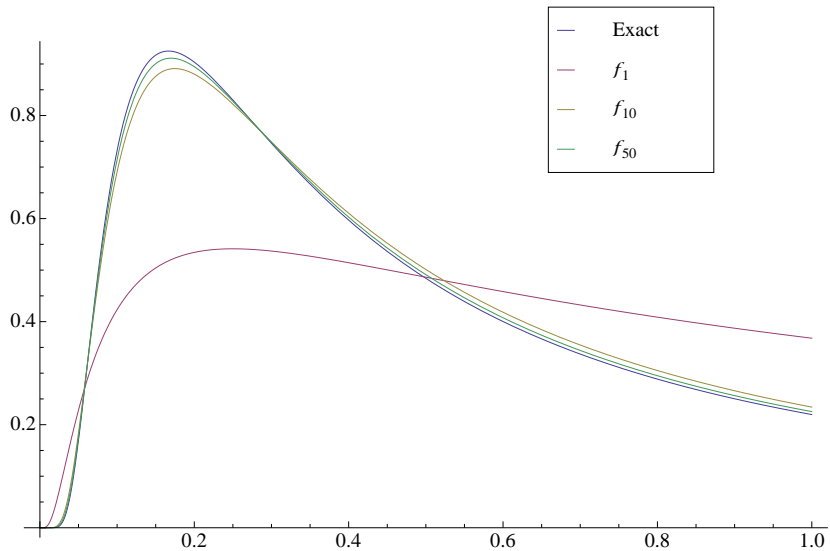}
\caption{ The exact form of $f$ together with the PW approximations $f_1$, $f_{10}$ and $f_{50}$ } \label{f:PWapprox}
\end{figure}

The convergence of $g_k(x)$ to $g(x)$ is slow in general.  To
illustrate this point, let $g$ be a PDF $f$ of an inverse Gaussian
distribution for which the Laplace transform is given by
$\tilde{f}(\lambda) = \exp(- \lambda^{-1/2})$ (in this case the
density $f(x)$ is known in closed form, see Section \ref{s:testit}).  In
Figure \ref{f:PWapprox}, we plot the exact formula for $f$ as well as
$f_k$ for $k=1,10 \ \& \ 50$.  Notice that even with $49$ derivatives
of $\tilde{f}$, the approximation is still poor.  In fact, this is a
general feature of the PW formula, as it has been shown (\cite{jagerman:1982}) that the errors $\epsilon_k(f;x) \equiv f(x) - f_k(x)$ have a power series expansion:
\begin{equation}\label{e:PWerror}
\epsilon_k(f;x) = \sum_{m=1}^\infty \frac{a_m(x)}{k^m}
\end{equation}
where the coefficients $a_m(m)$ are given by
\begin{equation*}
a_m(x) = \sum_{j=1}^m f^{(j+m)}(x) x^{j + m} \frac{ d(j + m,m)}{(j + m)!},
\end{equation*}
where $d(j+m,m)$ are the associated Stirling numbers of the first kind
(\cite{riordan:2002}, Chapter 4, section 4).
Thus, the convergence of $f_k$ is in general $O(1/k)$ corresponding to
the $m=1$ term.  This method alone is thus inadequate for computing the inverse Laplace transform to a high level of precision.

To obtain high precision, it is necessary to couple the PW approximations with a convergence acceleration method which extrapolates the limit in (\ref{e:PWformula}) based on a finite collection of terms in the sequence.  Thus, our method for computing the PDF $f$ at $x>0$ involves two components:

\begin{itemize}

\item[1.] For some sequence $k_1 < k_2 < \dots < k_N$, compute the approximates $f_{k_i}(x)$.

\item[2.]  Use the points $(k_1,f_{k_1}(x)), \dots, (k_N,f_{k_N}(x))$ to extrapolate $\lim_{k\rightarrow \infty} f_{k}(x)$.

\end{itemize}
And for the CDF, replace $f$ with $F$ above.  We focus first on step 1, which
in sight of (\ref{e:fk}) involves computing high derivatives of $\psi$ or
$\Psi$.   Next, we'll discuss two methods of convergence acceleration
to address step 2.

\medskip

\noindent \textbf{Remark:}  Our method produces the best results if the PDF and/or
CDF are smooth functions for $x > 0$.  That being said,  this method
still produces useful results in the non-smooth case,
however it may fail to converge near points where the
function lacks smoothness.  This method should {\it not} be used for
distributions which contain atoms, such as the Poisson distribution,
or even distributions whose density lacks smoothness, such as a compound
Poisson distributions with bounded jump distribution.  For such
distributions another method based on
numerically solving a Kolmogorov-Feller forward equation, is more
suitable, see \cite{veillette:2009a}.

\medskip

\subsection{Computing derivatives}

\subsubsection*{Derivatives of $\psi$}  We begin by presenting
two methods of computing derivatives of $\psi(\lambda) = \exp(-\phi(\lambda))$.  The first method is
based on the simple observation that
\begin{equation*}
 \psi'(\lambda) = - \phi'(\lambda) \psi(\lambda).
\end{equation*}
Thus, Leibnitz's formula implies that for any $k \geq 1$, the
$(k+1)^{th}$ derivative of $\psi$ is given by
\begin{equation}\label{e:derofpsi}
\psi^{(k+1)}(\lambda) = -  \sum_{j = 0} ^ k { k \choose j} \psi^{(j)}(\lambda) \phi^{(k + 1 - j)}(\lambda).
\end{equation}
Notice that for $\lambda$ fixed, the values of $\psi(\lambda),\psi'(\lambda),\dots,\psi^{(k)}$ can be
computed recursively using (\ref{e:derofpsi}) if one has first computed $\phi^{(j)}(\lambda)$ for
$j = 1,2,\dots,k+1$.  

Alternatively, the $k^{th}$ derivative of $\psi(\lambda)
= \exp(-\phi(\lambda))$ can be given directly in terms of the
derivatives of $\phi$ using Fa\'{a} di Bruno's formula
(\cite{riordan:2002}, Chapter 2, section 8),
which yields in our case
\begin{equation*}
\psi^{(k)}(\lambda) = \frac{d^k}{d \lambda^k} e^{-\phi(\lambda)} = e^{-\phi(\lambda)}  \sum \frac{k
  !}{m_1! m_2! \dots m_k!} \left(
  \frac{-\phi^{(1)}(\lambda)}{1!} \right)^{m_1}  \dots
\left( \frac{-\phi^{(k)}(\lambda)}{k!} \right) ^{m_k} ,
\end{equation*}
where the  sum runs over all $m_1,\dots,m_k$ such that $m_1 + 2 m_2
+ \dots + k m_k = k$.  This can expressed more simply in terms of the $k^{th}$ {\it complete Bell
  polynomial,} $B_k(x_1,x_2,\dots,x_k)$, as
\begin{equation} \label{e:belleqn}
\psi^{(k)}(\lambda) = e^{-\phi(\lambda)}
B_k(-\phi^{(1)}(\lambda),-\phi^{(2)}(\lambda),\dots,-\phi^{(k)}(\lambda)
).
\end{equation}
See \cite{pitman:2002} or \cite{roman:1984}.  Furthermore, $B_k$ can be given in terms of the following $k \times k$
determinant: 
\begin{equation}\label{e:belldet}
B_k(x_1,\dots,x_k) = \left| \begin{array}{cccccc}
x_1 & {k- 1 \choose 1} x_2 &  {k- 1 \choose 2} x_3 &   {k-
  1 \choose 3} x_4 & \dots &
x_k\\
-1 & x_1 & {k- 2 \choose 1} x_2 & {k- 2 \choose 2} x_3 &
\dots &  x_{k-1}  \\
0 & -1 & x_1 & {k- 3 \choose 1} x_2 & \dots & x_{k- 2}  \\
0 & 0 & -1 & x_1 & \dots & x_{k- 3} \\
\vdots & \vdots & \vdots  & \vdots & \ddots& \vdots  \\
0  & 0 & 0 & \dots & -1& x_1
\end{array} \right|.
\end{equation}
This provides another way for computing $\psi^{(k)}$.

  Both methods above require the derivatives of the Laplace exponent $\phi$.  Fortunately,  these can be computed directly, since the LK formula implies that for any $\lambda > 0$,
\begin{equation}\label{e:derofphi}
\phi^{(n)}(\lambda)  = \begin{cases} \displaystyle \int_0^\infty (1-e^{-\lambda u}) \Pi(du), \quad &n = 0 \\ \\
 \displaystyle  (-1)^{n+1} \int_0^\infty u^n e^{-\lambda u} \Pi( du), \quad &n \geq 1  \end{cases}
\end{equation}
Notice that passing derivatives through the integrals above is justified by the integrability assumption (\ref{e:intassum}), since for any $n \geq 1$,
$x^n e^{-\lambda x} \leq C_{n,\lambda} (1 \wedge x) $  uniformly on
$\{x>0\}$ for some constant $C_{n,\lambda}>0$.   For many examples,
these integrals have a closed form expression and can be computed
easily.  If a closed form is not available, many numerical integration
methods are effective as these integrands are non-negative,
exponentially decaying functions.  In this case, it is important to
use a small  {\it relative} error tolerance, as $\phi^{(n)}(\lambda)$
can become extremely small for large $\lambda$.

\subsubsection*{Derivatives of $\Psi$}  For $\Psi$ defined in (\ref{e:derPSI}), we again apply Leibnitz's formula to (\ref{e:derPSI}) and obtain
\begin{equation*}
\Psi^{(k)}(\lambda) = \tilde{F}^{(k)}(\lambda) = \frac{d^k}{d\lambda^k} \frac{\psi(\lambda)}{\lambda}= \sum_{j = 0}^{k} { k \choose j}  (-1)^j \frac{j!}{\lambda^{j+1}} \psi^{(k-j)}(\lambda),
\end{equation*}
 When
combined with the Post-Widder formula (\ref{e:fk}), we see a good deal of cancellation:
\begin{align}
F_k(x) &= \frac{(-1)^k}{(k-1)!} \left( \frac{k}{x} \right)^k
\tilde{F}^{(k-1)}(k/x) \nonumber \\
&= \frac{(-1)^{k-1}}{(k-1)!} \left( \frac{k}{x} \right)^{k} \Psi^{(k-1)}(k/x) \nonumber\\ &=  \frac{(-1)^{k-1}}{(k-1)!} \left( \frac{k}{x} \right)^{k} \left[  \sum_{j = 0}^{k-1} { k-1 \choose j}  (-1)^j \frac{j!}{(k/x)^{j+1}} \psi^{(k-1-j)}(k/x) \right] \nonumber \\
&= \sum_{j = 0}^{k-1} \frac{ (-1)^{k+j-1}}{(k - j - 1)!} \left( \frac{k}{x} \right)^{k-j-1} \psi^{(k-j-1)}(k/x) \label{e:PWCDF}
\end{align}
Note that this requires the values
$\psi(\lambda),\psi'(\lambda),\dots,\psi^{(n)}(\lambda)$, which can be
computed using (\ref{e:derofpsi}) and (\ref{e:derofphi}) hence the recursive formula (\ref{e:derofpsi}) is more
convenient here than (\ref{e:belleqn}) since it computes all derivatives of $\psi$.

\subsubsection*{Other useful derivatives}

In some applications, the derivatives of the PDF are also of
interest.  Assuming the first $q-1$ derivatives of $f$ vanish at $0$, i.e. $f(0) = f'(0) = \dots =
f^{(q-1)}(0) = 0$, then $f^{(q)}(x)$ has Laplace
transform
\begin{eqnarray*}
\widetilde{f^{(q)}}(\lambda) &=& \lambda^q \psi(\lambda) - \lambda^{q-1}
f(0) - \lambda^{q-2} f'(0) - \dots - f^{q-1}(0) \\
&=& \lambda^q \psi(\lambda) ,
\end{eqnarray*}
where $\psi(\lambda) = \tilde{f}(\lambda)$.  In which case, Leibnitz's formula gives
\begin{equation}
\frac{d^k}{d \lambda^k}  \widetilde{f^{(q)}}(\lambda)  =\frac{d^k}{d \lambda^k} \lambda^q \psi(\lambda) = \sum_{j=0}^{q \wedge (k-1)}
{k \choose j} \frac{q!}{(q- j)!} \lambda^{q- j} \psi^{(k-
  j )} ( \lambda) 
\end{equation}
Thus, (\ref{e:derofpsi}) and (\ref{e:derofphi}) can be applied to compute the derivatives of $f$.

\subsection{Extrapolation}

We want to compute the PDF $f(x)$ or CDF $F(x)$ by approximating the
limit as $k \rightarrow \infty$ in the Post Widder formula (\ref{e:PWformula}).     By letting $h = k^{-1}$, (\ref{e:PWerror}) implies that $f(x)$ can be written as
\begin{equation}
f(x) = f_k(x) + \sum_{m=1}^\infty a_m h^{m}
\end{equation}
We shall truncate the error series $\sum_{m=1}^\infty$ at $m = N-1$
and consider non-consecutive $k_1 <
k_2 < \dots < k_n$.  If $h_i =
k_i^{-1}$, $i=1,\dots,N$, we can use the points
$(h_1,f_{k_1}(x)),(h_2,f_{k_2}(x)),\dots,(h_N,f_{k_N}(x))$ to estimate
the first $N-1$ unknown coefficients $a_m$.  Indeed, the (approximate
) $N \times N$ system
\begin{equation}\label{e:sysas}
 f_{k_i}(x)  \approx  f(x) - a_1 h_j - a_2 h_j^2 - \dots - a_n
 h_j^{N-1}, \quad j=1,2,\dots,N
\end{equation}
can be solved for the $N$ unknowns $f(x),a_1,\dots,a_{N-1}$.  Since we
are only interested in $f(x)$, it is enough to only compute the first
row of the inverse of the matrix corresponding to the system
(\ref{e:sysas}).  This is essentially done in the polynomial
extrapolation method described below in order to  obtain an accurate approximation of the limit $\lim_{k \rightarrow
  \infty} f_{k}(x)$.  There is a vast literature on such extrapolation
methods and the errors associated to them, see \cite{joyce:1971} for a
review.  

Below, we briefly review two techniques for approximating $f(x)$ also
used in \cite{frolov:1998}.  We
call these methods of extrapolation.  The first is based on polynomial interpolation,
and the second based on rational function (Pad\'{e}) interpolation.
Not surprisingly, the rational extrapolation provides faster
convergence in many cases, however we found it more susceptible to
numerical instability for larger $N$ and $k_N$.   Therefore, the
``better'' choice depends on the particular example and the desired
accuracy.

\subsubsection*{Polynomial extrapolation}
Given the points $(h_1,f_{k_1}(x)),(h_2,f_{k_2}(x)),\dots,(h_N,f_{k_N}(x))$, the $N-1$ degree Lagrange polynomial $P_N(x)$ which passes through these points is given by
\begin{equation*}
P_N(h) = \sum_{i=1}^N \left(  \prod_{j \neq i} \frac{h - h_j}{h_i - h_j} \right) f_{k_i}(x).
\end{equation*}
Thus,  taking $h = 0$ (or, $k=\infty$) in this polynomial provides the approximation 
\begin{equation*}
f(x) \approx P_N(0) = \sum_{i=1}^N c_i f_{k_i}(x), \quad \mbox{where} \quad c_i = \prod_{j\neq i} \frac{-h_j}{h_i - h_j} = \frac{ (-1)^{N-1} k_i^{N-1}}{\prod_{j \neq i} (k_j - k_i)}.
\end{equation*}
We thus use a linear combination of the approximations $f_{k_i}(x)$ to obtain
a more accurate approximation of $f(x)$.

\subsubsection*{Rational extrapolation}
As an alternative to polynomial extrapolation, on may instead fit a rational
function to the points $( h_1 , f_{k_1}(x) )$, $\dots$,$( h_N , f_N{ k_N
} (x))$ of the form
\begin{equation*}
R_N(h) =  \frac{P_{\mu}(h)}{Q_{\nu}(h)},
\end{equation*}
where $P_{\mu}(h) = \sum_{i=1}^\mu p_i h^i$ is a polynomial of degree $\mu =
\lfloor N/2 \rfloor$ and $Q_\nu(h) = \sum_{i=0}^\nu q_i h^i$ is a polynomial of degree
$\nu = N - \lfloor N/2 \rfloor$.   Here,  the coefficients $p_i,q_i$ are
chosen so $R_N(h_i) = f_{k_i}(x)$.   When implementing this method, these
coefficients are not computed directly like in the polynomial extrapolation
case, but instead $R_N(0)$ is computed iteratively in a triangular array using the following
recursive formulas (see \cite{joyce:1971}, section 13)
\begin{align*}
& R_{-1}^i = 0,\qquad  R_{0}^i = f_{k_i}(h_i), \quad  i \geq 1 \\ 
& R_{m}^i = R_{m-1}^{i+1} + \frac{ R^{i+1}_{m-1} - R^i_{m-1}}{h_i \Theta^i_m -
    h_{i+m}}, \quad  \mbox{where} \quad \Theta_m^i = \frac{R_{m-1}^i -
    R_{m-1}^{i+1}}{R^{i+1}_{m-1} - R^{i+1}_{m-2}}.
\end{align*}  
We then take $f(x) \approx R^1_N = R_N(0)$ as our approximation since
$h = 1/k \rightarrow 0$ corresponds to $k\rightarrow \infty$. Notice that unlike the
polynomial extrapolation, this transformation is {\it nonlinear} in the $f_{k_i}(x)$.

\subsubsection*{Error bounds}

One major drawback to  using extrapolation techniques is that the error made in
the approximation is difficult to bound analytically (see, for instance,
Theorem 3 in \cite{bulirsch:1966} and equation (13) in \cite{frolov:1998}).
One can however, obtain an asymptotic error bound for a rational or polynomial
approximation.  Let $h_1 < h_2 < h_3 \dots$ be a sequence for which the approximant
$f_{k_i}(x)$ can be computed with $k_i = 1/h_i$.  Also, let $P_N(0)$  be the approximation of $f(x)$ obtained by using polynomial
 extrapolation with $h_1,h_2,\dots,h_N$, with $N \geq 2$ (all the following
 formulas also hold with rational extrapolation by replacing $P_N(0)$ with $R_N(0)$).  In
\cite{bulirsch:1966}, Bulirsch and Stoer (BS) construct a {\it second} estimate
of $f(x)$,  $\widetilde{P}_N(0)$, with the property that 
\begin{equation}\label{e:asyerrorbound}    
\lim_{N \rightarrow \infty} \frac{ P_N(0) - f(x) }{\widetilde{P}_N(0) - f(x)} = -1.
\end{equation}
Thus, asymptotically, $\widetilde{P}_N(0)$ is as good an approximation to $f(x)$ as $P_N(0)$,
except that  it approaches $f(x)$ from the opposite direction.  This second
approximation is obtained using a linear combination of the form 
\begin{equation}\label{e:newp0}
\widetilde{P}_N(0)
= (1 + \alpha)
P_{N+1}(0) - \alpha P_N(0),
\end{equation} 
and choosing the constant $\alpha$ to change the
sign of the leading term in the error $P_N(0) - f(x)$. BS show
this constant $\alpha$
is given by
\begin{equation}\label{e:defofalpha} 
\alpha = 1 + \frac{2}{\frac{h_{1}}{h_{N+1}} - 1}.
\end{equation}
This allows us to construct a numerical bound on the relative error made in this
method.  From (\ref{e:asyerrorbound}), we have
\begin{equation}\label{e:goto2}
\frac{ P_N(0) -
  \widetilde{P}_N(0) }{ f(x) - \widetilde{P}_N(0)} =  \frac{ P_N(0) - f(x) }{f(x) - \widetilde{P}_N(0)}  +1  \longrightarrow
2.
\end{equation}
Fix any $\eta \in (\frac{1}{2},1)$.  From (\ref{e:goto2}) and since $\eta>\frac{1}{2}$, there exists $N_0(\eta)$ large enough
such that for $N>N_0(\eta)$,  $|f(x) -
\widetilde{P}_N(0)| < \eta |P_N(0) -
  \widetilde{P}_N(0)|$.  Using a similar argument, there also exists $M_0(\eta)$
  large enough such that for $N>M_0(\eta)$, $|P_N(0) - f(x)| < \eta |P_N(0) -
  \widetilde{P}_N(0)|$.   Moreover, since $P_N(0) \rightarrow
  f(x)$ and $\eta < 1$, we can pick $N_1(\eta)$ large enough such that $|P_N(0)/f(x)|
  < 1/\eta$ for all $N>N_1(\eta)$.  

We thus can refine our estimate as $f(x) \approx (P_N(0) +
\widetilde{P}_N(0))/2$.  If $N>\max(N_0(\eta),M_0(\eta))$,
\begin{eqnarray*}
\left| \frac{P_N(0) + \widetilde{P}_N(0)}{2} - f(x)  \right| \leq
\frac{ | P_N(0)-f(x)| }{2  } + \frac{ |f(x) -\widetilde{P}_N(0) |
}{2 } \leq \eta |P_N(0) - \tilde{P}_N(0) | \leq  |P_N(0) - \tilde{P}_N(0) |,
\end{eqnarray*}
and if $N>\max(N_0(\eta),M_0(\eta),N_1(\eta))$,
\begin{eqnarray*}
\frac{\left| \frac{P_N(0) + \widetilde{P}_N(0)}{2} - f(x)  \right|}{f(x)} &\leq&
\frac{ | P_N(0)-f(x)| }{2 f(x) } + \frac{ |f(x) -\widetilde{P}_N(0) |
}{2 f(x)}  \\
&\leq&   \eta  \left(\frac{P_N(0)}{f(x)} \right)\frac{  | P_N(0) - \tilde{P}_N(0) |  }{P_N(0)} \\
&\leq&   \frac{ | P_N(0) -  \tilde{P}_N(0) | }{P_N(0)} .
\end{eqnarray*}

Therefore, 

\begin{thm}
Suppose $f(x) \neq 0$ and let $P_N(0)$ be an approximation of $f(x)$ obtained using either
polynomial extrapolation or rational extrapolation, and let
$\widetilde{P}_N(0)$ be defined as in (\ref{e:newp0}).  For $N$ large
enough, we have 
\begin{equation}\label{ebound}
 \left| \frac{P_N(0) + \widetilde{P}_N(0)}{2} - f(x)  \right|
  \leq   | P_N(0) -
\tilde{P}_N(0) | \quad \mbox{and} \quad  
\frac{ \left| \frac{P_N(0) + \widetilde{P}_N(0)}{2} - f(x)  \right|
}{f(x) }  \leq   \frac{| P_N(0) -
\tilde{P}_N(0) | }{P_N(0)},
\end{equation}
where the right hand sides converge to 0 as $N \rightarrow \infty$.
\end{thm}

This provides a numerical bound on the absolute and  relative errors of the approximation and suggests increasing $N$
until either $|P_N(0)-\tilde P_N(0)|$  or $|P_N(0)-\tilde P_N(0)| / P_N(0) $ is smaller than a prescribed value.

\section{Implementation}\label{s:implementation}
We begin with a description of the algorithm for computing the
PDF and CDF based using the recursive formula (\ref{e:derofpsi})
and the polynomial extrapolation method.  This method is implemented
in MATLAB\footnote{MATLAB is computational software package developed by The
  Mathworks} and {\it Mathematica}\footnote{Mathematica is a
  computational software package developed by Wolfram Research}.  The method using rational extrapolation is
also an option in the MATLAB implementation, but is not described
here.  

\medskip

\medskip

\noindent  \fbox{ {\bf Algorithm for computing the PDF $f$ at $x>0$} }

\medskip

\begin{itemize}

\item[{\bf 1.}] Choose the sequence $1 = k_0 < k_1 < k_2 < \dots < k_N$ and a relative error tolerance $\epsilon > 0$.  We found $k_j= 10 j$, with $N$ lying around $8$ to be effective for $\epsilon = 10^{-6}$.  Initialize the following arrays:
\begin{equation*}
\begin{tabular}{ccc}
\mbox{Variable} & \mbox{Size} & \mbox{Purpose} \\
$D_{\phi}$ & $(k_N-1) \times N$ & \mbox{Holds derivatives of $\phi$} \\
$D_{\psi}$ & $k_N \times N$ & \mbox{Holds derivatives of $\psi$} \\
$P$ & $1 \times N$ & \mbox{Holds Post-Widder approximations of $f$} \\
$f$ & $1 \times N$ & \mbox{Holds extrapolated approximations of $f$}
\end{tabular}
\end{equation*}
Set $j=1$.

\item[{\bf 2.}]  Compute $(D_{\psi})_{1,j}$ as
\begin{equation*}
(D_\psi)_{1,j} = \psi(k_j/x)
\end{equation*}

\item[{\bf 3.}] For $i=k_{j-1},\dots,k_j-1$, compute the $i^{th}$ row of $D_\phi$ as
\begin{equation*}
(D_\phi)_{i,j} = \begin{cases}  \phi^{(i)}(k_j/x) \quad &i \leq k_j \\
0 \quad &\mbox{otherwise} \end{cases}, \quad j=1,2,\dots,N.
\end{equation*}
$\phi^{(i)}$ is computed using the integrals in (\ref{e:derofphi}).  The $(i+1)^{th}$ row of $D_\psi$ is then computed as
\begin{equation*}
(D_\psi)_{i+1,j} = \begin{cases} \displaystyle   \sum_{r = 1}^{i} {i \choose r} (D_\psi)_{r,j} (D_\phi)_{i + 1 - r,j}   \quad & i \leq k_j\\
0 \quad &\mbox{otherwise} \end{cases}, \quad j=1,2,\dots,N.
\end{equation*}

\item[{\bf  4.}] Compute the $j^{th}$ Post-Widder approximation:
\begin{equation*}
P_j = \frac{(-1)^{k_j - 1}}{(k_j - 1)!} \left(\frac{k_j}{x}\right)^{k_j} (D_{\psi})_{k_j,j} .
\end{equation*}

\item[{\bf  5.}] If $j=1$, set $f_1 = P_1$, $j = 2$, and go to step 3.  Otherwise, compute the $j^{th}$ extrapolation as
\begin{equation*}
f_j = \sum_{r = 1}^j c_r P_r, \quad \mbox{where} \quad c_r = \prod_{\ell= 1,\ell \neq r}^r \frac{k_r}{k_\ell - k_r}.
\end{equation*} 
Set $\alpha = 1 + 2(k_j/k_1 + 1)^{-1}$.  If $|(1+\alpha) (f_j -f_{j-1})|/f_{j} < \epsilon$ or if $j = N$, return $f_j$.  Otherwise, set $j = j+1$ and go to step 3.

\end{itemize}

The method for computing the CDF is similar, except that (\ref{e:PWCDF}) is used instead for step 4.

\medskip

\medskip

\noindent  \fbox{ {\bf Algorithm for computing the CDF $F$ at $x>0$} }

\medskip

\noindent Follow steps for PDF computation above except replace step 4. with
\begin{itemize}

\item[{\bf 4.$^\star$}] Compute the $j^{th}$ Post-Widder approximation:
\begin{equation*}
P_j  = \sum_{r = 0}^{k_j-1} \frac{ (-1)^{k_j+r-1}}{(k_j - r - 1)!} \left( \frac{k_j}{x} \right)^{k_j-r-1} (D_\psi)_{k_j-r,j}
\end{equation*}

We'll conclude this section with a collection of remarks regarding implementing this procedure.  

\subsubsection*{Remarks}

\begin{itemize}

\item[1.] Binomial coefficients are reused multiple times in step 3.  We found it useful to compute and store rows of Pascal's triangle as needed.  

\item[2.]  In MATLAB, all computations above can be ``vectorized'' to
  maximize speed. 

\item[3.]  It is essential to compute the derivatives of $\phi$ to as many {\it significant} digits of accuracy in step 3 as you want in your final result.  In many examples, a closed form for $\phi^{(n)}$ can be found.  See Section \ref{s:examples} for examples where these integrals are computed without a closed form.   

\item[4.]  To avoid overflow, ratios of large numbers like those seen in steps 4.~and 4.$^\star$~should be computed to incorporate reduction.  For example, using identities such as
\begin{equation*}
\frac{1}{(k_j -1)!} \left( \frac{k_j}{x} \right)^{k_j} = \exp\left( k_j \log(k_j/x) - \sum_{j=1}^{k_j-1} \log(j) \right)
\end{equation*}
will alleviate overflow.  

\item[5.] In double precision arithmetic, it is often impossible to
  take $k_N$ much bigger than $100$ as underflow, overflow and
  numerical instability become unavoidable.  In most cases, adequate
  convergence is met much before such large $k$ values are needed,
  however if one must go further,  it will become necessary to work
  with a high-precision arithmetic package.  {\it Mathematica}, for example, has
  powerful multiprecision capabilities. If one requires more speed, a
  well developed and documented library for multiple precision
  arithmetic for C++/Fortran is available at \url{http://crd.lbl.gov/~dhbailey/mpdist.} 

\item[6.]  When the true value of $f$ is very close to $0$ ($f(x) \ll
  10^{-8}$) the extrapolation procedure might ``overshoot'' 0 and
  return a negative value for the density.  In this case,  $0$ is a
  good approximation as you can be sure that the PDF or CDF takes an
  extremely small value in this case.  And similarly when the CDF takes a value greater than 1.

\end{itemize}

\end{itemize}

\section{Testing the method}\label{s:testit}

In this section we use {\it Mathematica} and consider two examples for which it is
possible to compute the PDF and CDF exactly.  The first example is the
chi-squared distribution $\chi^2$ with one degree of freedom.  For this distribution,
\begin{align*}
&f_{\chi^2}(x) = \frac{1}{\sqrt{2 \pi}} x^{-1/2} e^{-x/2}, \quad F_{\chi^2}(x) = \mbox{erf}\left(\sqrt{\frac{x}{2}}\right) \\
&\psi_{\chi^2}(\lambda) = \tilde{f}_{\chi^2}(\lambda) =  (1 + 2 \lambda)^{-1/2} \\
&\phi_{\chi^2}(\lambda) = -\log \psi_{\chi^2}(\lambda)= \frac{1}{2} \log(1 + 2 \lambda) =  \int_0^\infty (1 - e^{-\lambda u}) \left[ \frac{e^{-u/2}}{2 u} \right] du. \\
& \phi^{(n)}_{\chi^2}(\lambda) =  (-1)^{n+1} \int_0^\infty u^n
e^{-\lambda n} \left[ \frac{e^{-u/2}}{2 u} \right] du  \\
&\qquad =   \frac{(-1)^{n+1}}{2} \int_0^\infty  u^{n-1} e^{-u(\lambda
  + 1/2)}  du = \frac{(-1)^{n+1} (n-1)!  }{2} \left( \frac{1}{2} + \lambda \right)^{-n}, \quad n \geq 1.
\end{align*}

\begin{figure}[ht]
\centering
\subfigure[Relative error for PDF of $\chi^2$ distribution]{
\includegraphics[scale=.5]{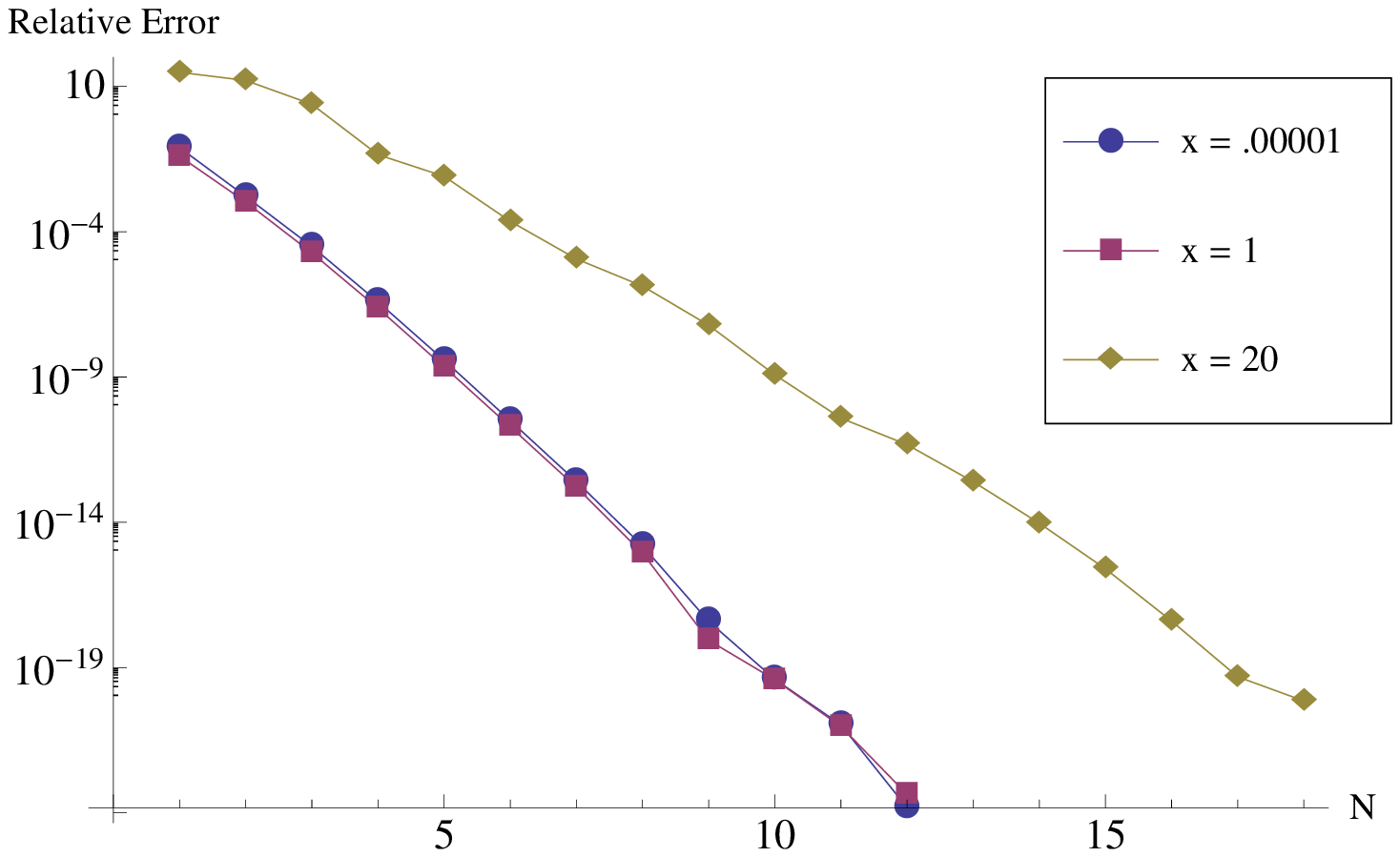}
}
\subfigure[Relative error for  CDF of $\chi^2$ distribution]{
\includegraphics[scale=.5]{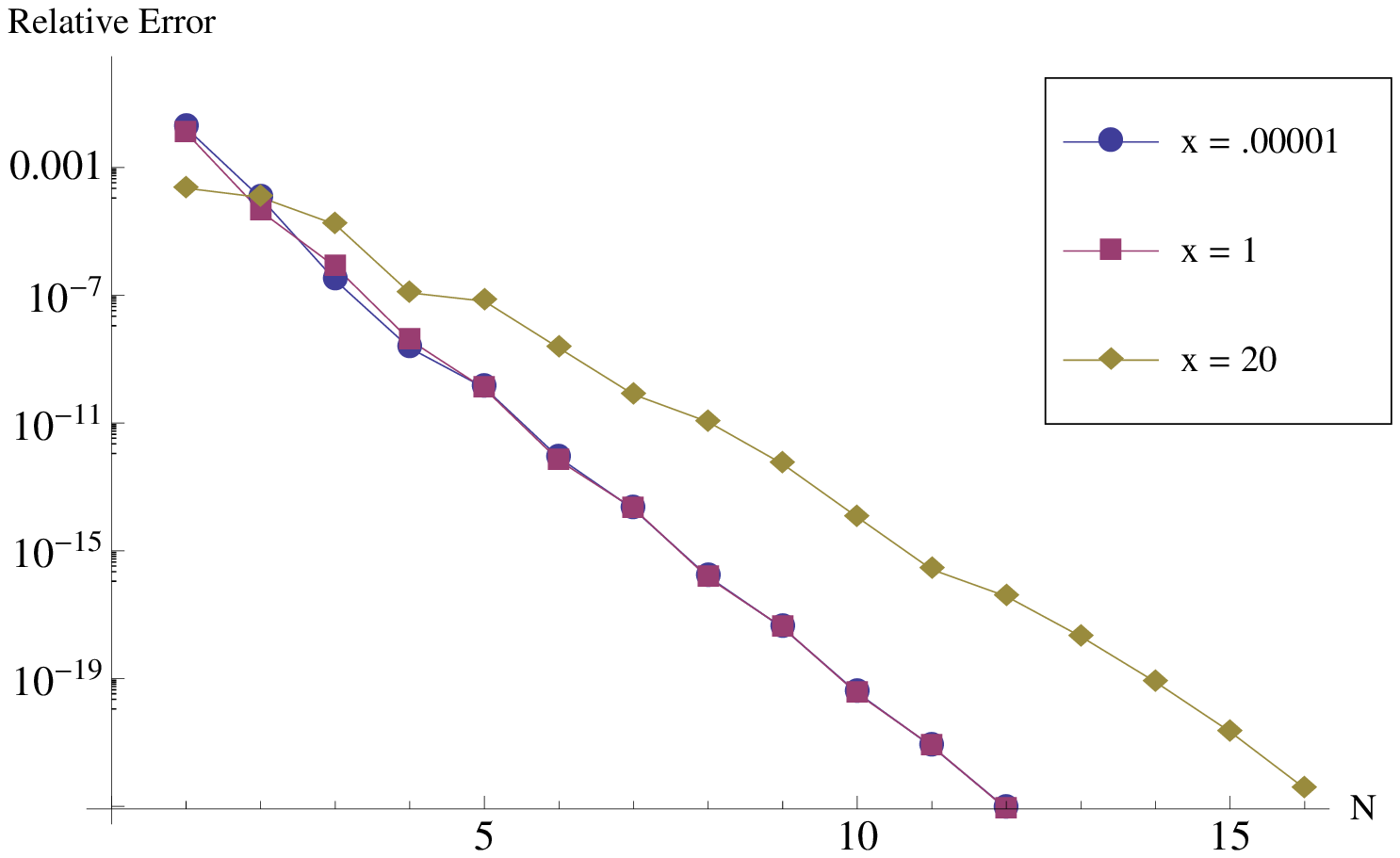}
}
\caption[Relative error of method with polynomial extrapolation]{ Plot
  of relative error of our method versus the truncation level $N$ for the PDF and CDF of
  the $\chi^2$ distribution at three values of $x$.  \label{f:C2errplots} }
\end{figure}

\begin{center}
\begin{table}
\begin{tabular}{|l|l|l|l|l|l|l|l|} 
\hline 
\multicolumn{7}{|c|}{Chi-squared} \\
\hline
$x$ & $10^{-5}$ & $10^{-1}$ & $1$ & $10$ & $20$ & $50$ \\ \hline \hline
\multirow{1}{*}{$f_{\chi^2}(x)$} &
 $1.26 \times 10^2$ &  1.20 &  0.242 & $8.50 \times 10^{-4}$ & $4.05 
 \times 10^{-6}$ &  $7.83 \times 10^{-13}$ \\ \hline
\multirow{1}{*}{$N_6$} &
 4 &  4 &  4 & 6 & 9 & 14 \\  \hline
\multirow{1}{*}{$N_{15}$} &
8 &8 & 9 & 12 & 15 & 21 \\
\hline \hline
\multirow{1}{*}{$F_{\chi^2}(x)$}  &
$2.52 \times 10^{-3}$ & 0.248 & 0.683 & 0.998 & $1 - 7.7 \times
10^{-6}$ & $1 - 1.54 \times
10^{-13}$ \\ \hline
\multirow{1}{*}{$N_{6}$} &
3 & 3 & 3 & 4 & 4 & 1 \\ \hline
\multirow{1}{*}{$N_{15}$} &
8 & 8 & 8 & 11 & 11 & 12 \\ \hline
\hline 
\multicolumn{7}{|c|}{Inverse Gaussian} \\
\hline
$x$ & $0.01$ & $0.02$ & $0.1$ & $1$ & $100$ & $1000$ \\ \hline \hline
\multirow{1}{*}{$f_{IG}(x)$} &
 $3.92 \times 10^{-9}$ &  $3.71 \times 10^{-4}$ &  0.732 & 0.220 & $2.81 \times 10^{-4}$ & $8.92 
 \times 10^{-6}$ \\ \hline
\multirow{1}{*}{$N_6$} &
 12 &  8 &  6 & 5 & 6 & 6 \\  \hline
\multirow{1}{*}{$N_{15}$} &
21 & 17 & 14 & 11 & 10 & 11 \\
\hline \hline
\multirow{1}{*}{$F_{IG}(x)$}  &
$1.53 \times 10^{-12}$ & $5.733 \times 10^{-7}$ & 0.025 & 0.480 &
0.944 & 0.982 \\ \hline
\multirow{1}{*}{$N_{6}$} &
12 & 9 & 6 & 4 & 4 & 3 \\ \hline
\multirow{1}{*}{$N_{15}$} &
21 & 17 & 13 & 10 & 8 & 8 \\ \hline
\end{tabular}
\caption{ Using $k_i = 10i$, this table shows the value of $N$
  required to obtain $6$ digits of accuracy ($N_6$) and 15 digits of
  accuracy ($N_{15}$) at the given values of $x$ using polynomial
  interpolation.  Notice that it becomes harder to approximate these
  functions to a small relative error when they become very small. }\label{t:Nval} 
\end{table}
\end{center}

\begin{figure}[ht]
\centering
\subfigure[Relative error for  PDF of IG distribution.]{
\includegraphics[scale=.5]{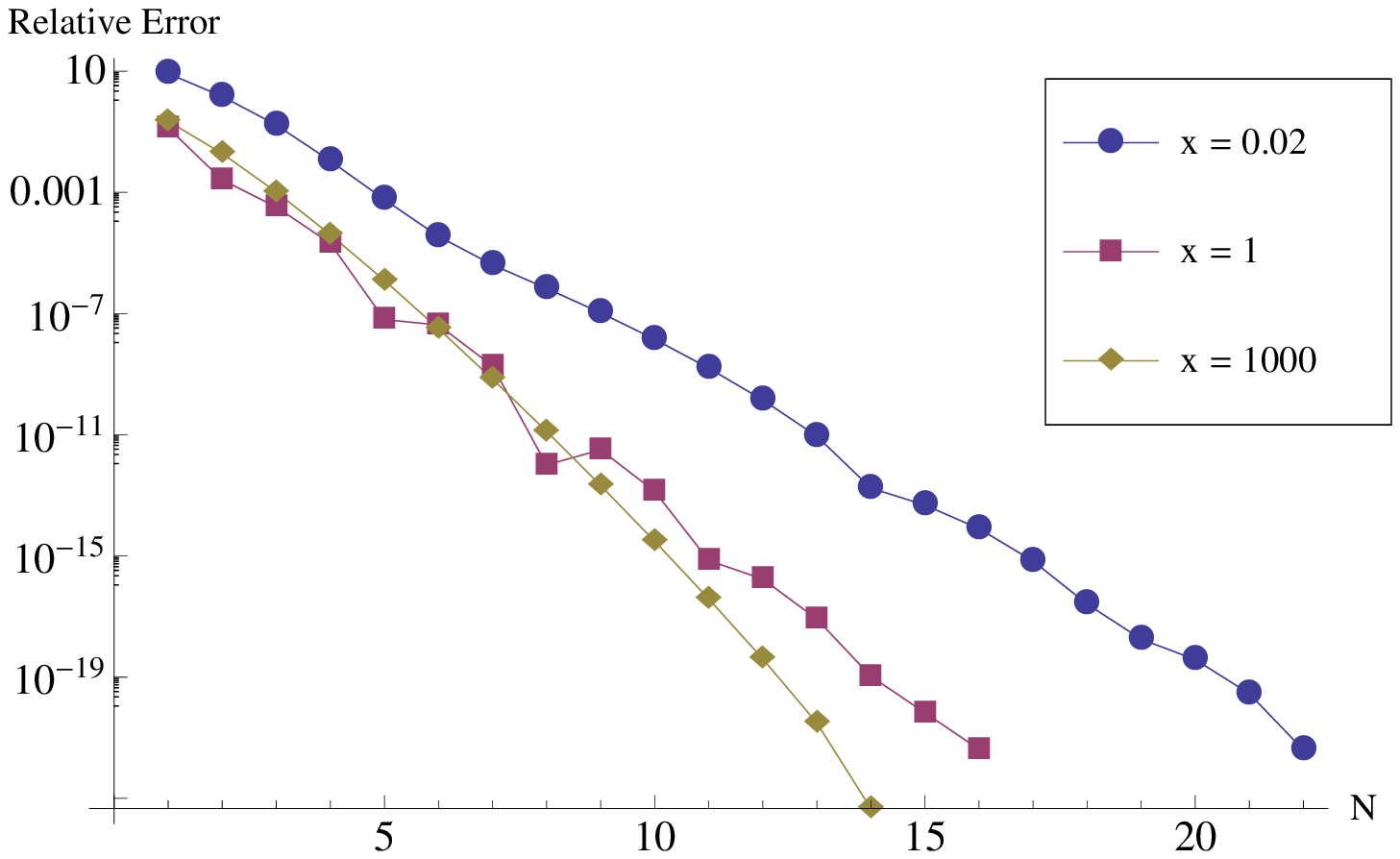}
}
\subfigure[Relative error for  CDF of IG distribution]{
\includegraphics[scale=.5]{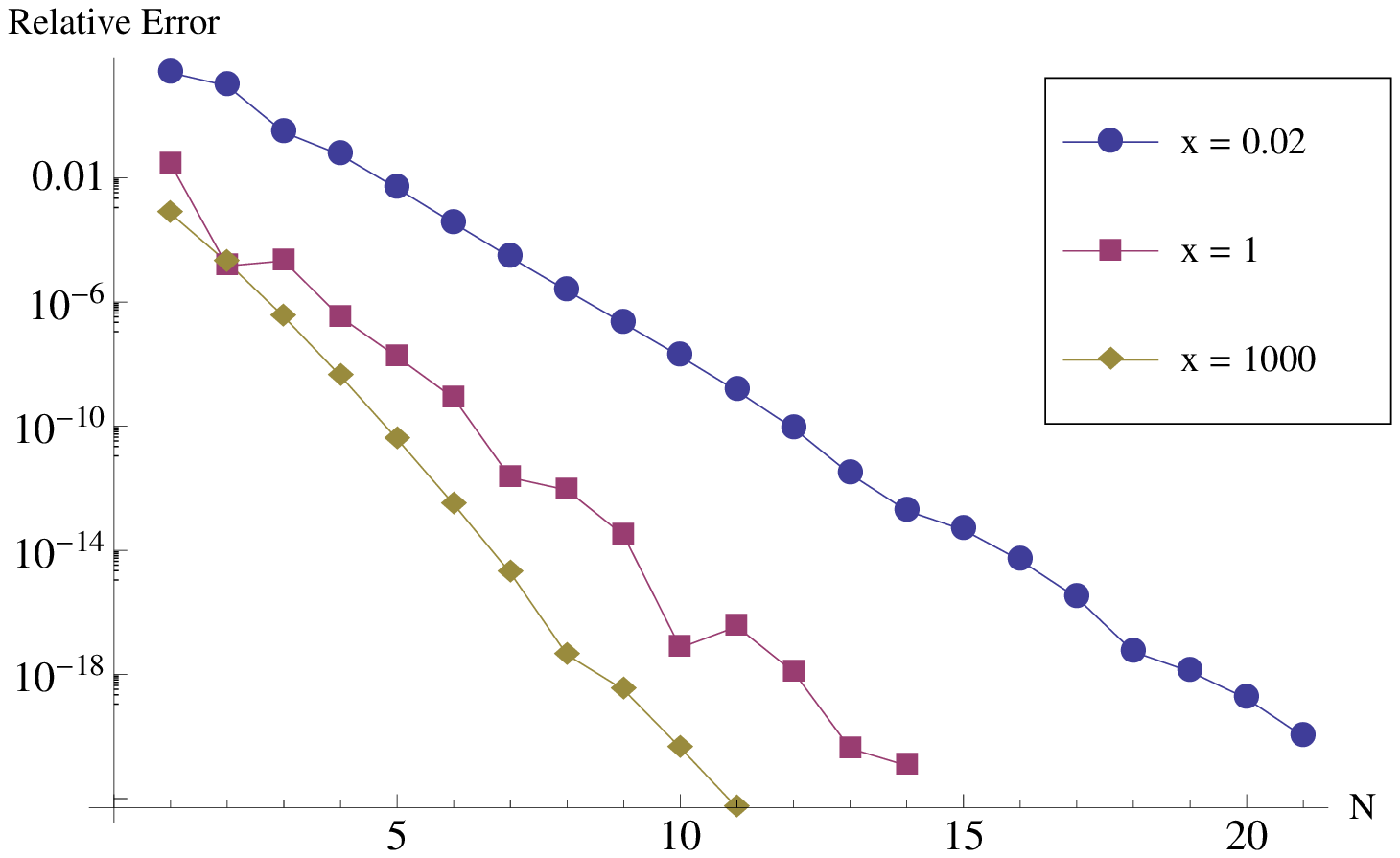}
}
\caption[Relative error of method with polynomial extrapolation]{ Plot
  of relative error of our method verses the truncation level $N$ for the PDF and CDF of
  the IG distribution at three values of $x$.  \label{f:IGeplots} }
\end{figure}

Results over a range of inputs is shown in Table \ref{t:Nval}.  We
show the $x$-values considered as well as the value of $N$ required to
obtain a relative precision of $10^{-6}$ and $10^{-15}$.  In Figure
\ref{f:C2errplots}, we also plot the relative error as a function of $N$.
Notice that the value of $N$ required is substantially higher when the
value of the PDF or CDF takes on very small values, in particular, when
$f_{\chi^2}(50) \approx 7.83 \times 10^{-13}$.  This is because
$f_k(x)$ over approximates $f$ at these points.

The second example we considered is the inverse Gaussian (IG)
distribution (\cite{Applebaum:2004}, Example 1.3.21), for which
\begin{align*}
& f_{IG}(x) = \frac{1}{\sqrt{4 \pi x^3}} e^{-1/(4x)}, \quad F_{IG}(x) = 1 - \mbox{erf}\left(\frac{1}{\sqrt{2 x}} \right) \\
& \psi_{IG}(\lambda) = \tilde{f}_{IG}(\lambda) =  e^{-\sqrt{\lambda}} \\
&\phi_{IG}(\lambda) = - \log \psi_{IG}(\lambda) =  \sqrt{\lambda} = \int_0^\infty (1 - e^{-\lambda u}) \left[ \frac{u^{-3/2}}{2\sqrt{\pi}} \right] du \\
&\phi^{(n)}_{IG}(\lambda) =\frac{ (-1)^{n+1}}{2 \sqrt{\pi}}
\int_0^\infty u^{n - 3/2} e^{-\lambda u} du = \frac{(-1)^{n+1}}{2 \sqrt{\pi}} \Gamma(n - 1/2)
\lambda^{1/2 - n} =  \frac{(-1)^{n+1}
  (2n-3)!  }{2^{2(n-1)}  (n-2)!  }  \lambda^{1/2 - n}   
\end{align*}
 Table \ref{t:Nval} also shows the values of $x$
considered, the value of $f_{IG}(x)$, and the value of $N$ needed to
obtain relative errors of $10^{-6}$ and $10^{-15}$. Figure \ref{f:IGeplots} shows how the relative error behaves as a function
of $N$ for three values of $x$.  Similarly to the $\chi^2$ distribution,
 $N$ is
largest when the value of the PDF is very close to $0$.

\section{Applications to ID distributions which are not known in closed form} \label{s:examples}

We will now consider examples of non-negative ID
distributions for which the PDF and CDF are not known in closed form.
We will compute them numerically using our method.  The resulting
plots of the PDFs and CDFs are obtained using MATLAB.  Each plot was
generated in about 1 second.  In order to apply our
method, we much first write the Laplace transform in LK form, and then
find $\phi(\lambda)$ and the derivatives $\phi^{(n)}(\lambda)$ by
computing the integrals in (\ref{e:derofphi}).  In almost every case
below, these integrals can be given in ``closed form'', which is to
say they can at least be written in terms of special functions for which there
are efficient methods for computation.  We will assume through out
that the drift $a = 0$ in (\ref{e:LK2}).

The following special functions will appear throughout this section,
we provide their definitions here for convenience:
\begin{align*}
& \mbox{Gamma function} & \quad & \Gamma(a) = \int_0^\infty z^{a-1}e^{- z}
dz, \quad a \neq 0,-1,-2,\dots \\
& \mbox{Lower incomplete gamma function}& \quad &\gamma(a,b) = \int_0^b
z^{a - 1} e^{-z} dz, \quad b \geq 0, \ a \neq 0,-1,-2,\dots \\
& \mbox{Upper incomplete gamma function}& \quad &\Gamma(a,b) =
\int_b^\infty z^{a-1} e^{-z} dz, \quad b>0, \ a \in \mathbb{R} \\
& \mbox{Entire exponential integral}& \quad &\mathrm{Ein}(a) =
\int_0^a \frac{1 - e^{-z}}{z} dz, \quad a \in \mathbb{R}  \\
& \mbox{Dilogatithm} & \quad &L_2(a) = \int_1^a \frac{\log(z)}{z - 1}
dz, \quad a > 0
\end{align*}
See \cite{oldham:2009}, chapters 25, 37, 43 and 45 for more discussion of
these functions and
methods for efficient computation.  For an integer $n \geq 1$,
we have $\Gamma(n) = (n-1)!$ and 
\begin{equation}\label{e:lincgam}
\gamma(n,b) = (n-1)! \left( 1 - e^{-b} \sum_{m=0}^{n-1} \frac{b^m}{m!} \right),  
\end{equation}
see \cite{gradshteyn:2007}, Equation
3.351.1. 

\subsection{Right-skewed stable distributions}

Here we consider the class of examples with Laplace exponent given by
\begin{equation}\label{e:st-lap}
\phi(\lambda) = \int_0^1 \lambda^\beta p(d \beta) = \int_0^\infty (1 -
e^{-\lambda u}) \left[ \int_0^1 \frac{u^{-\beta - 1}}{\Gamma(-\beta)}
  p(d \beta) \right] du,
\end{equation}
where $p$ is a measure supported on $(0,1)$.  The L\'{e}vy measure
$\Pi$ has a density is given by
\begin{equation*}
\Pi'(u) = \int_0^1  \frac{ u^{-\beta - 1}}{\Gamma(-\beta)} p(d\beta).
\end{equation*}
 Mixtures of this form are considered in
\cite{Meerschaert1:2004} and \cite{Meerschaert:2006}, in which these
distributions are used in models of anomalous diffusion.

The Laplace exponent is expressed in (\ref{e:st-lap}) and the derivatives $\phi^{(n)}$ for $n \geq 1$ can be computed by using (\ref{e:derofphi}), a change in the order of integration, and the definition of the gamma function:
\begin{equation*}
\phi^{(n)}(\lambda) = (-1)^{n+1} \int_0^\infty u^n e^{-\lambda u}
\left(\int_0^1 \frac{ u^{-\beta - 1}}{\Gamma(-\beta)} p(d\beta)\right)
du = (-1)^{n+1}
\int_0^1 \frac{\lambda^{\beta-n} \Gamma(n - \beta)}{\Gamma(-\beta)} p(d \beta)\quad n \geq 1.
\end{equation*}
Since $n$ takes only integer values, $\Gamma(n - \beta) = \Gamma(-\beta) \prod_{m=0}^{n-1} (m - \beta)$,  and so the above can be simplified further as
\begin{eqnarray}
\phi^{(n)}(\lambda) =  (-1)^{n+1} \int_0^1 \frac{\lambda^{\beta-n}
  \Gamma(n - \beta)}{\Gamma(-\beta)} p(d \beta) &=&  (-1)^{n+1} \int_0^1 \lambda^{\beta - n} \prod_{m=0}^{n-1}(m - \beta) p(d\beta)\nonumber \\
&=& -\frac{1}{\lambda^n} \sum_{m=1}^n S_n^{(m)} \int_0^1 \beta^m \lambda^{\beta} p(d \beta)\nonumber \\
&=& -\frac{1}{\lambda^n} \sum_{m=1}^n c_m(\lambda) S_n^{(m)} \label{e:asphin}
\end{eqnarray}
where $c_m(\lambda) = \int_0^1 \beta^m \lambda^{\beta} p(d\beta)$ and
$\{ S_n^{(m)} \}$, $n \geq 0$, $0 \leq m \leq n$ are the Stirling
numbers of the first kind (\cite{oldham:2009}, page 162).  These are
such that $\prod_{m=0}^{n-1} (\beta - m) = \sum_{m=1}^{n} S_n^{(m)}
\beta^m$,  and can computed with a triangular array similarly to Pascal's triangle using the recursion formula
\begin{align*}
& S_{0}^{(0)} = 1, \quad S_{n}^{(0)} = 0, \quad n \geq 1 \\
& S_n^{(m)} = S_{n-1}^{(m-1)} - (n-1)  S_{n-1}^{(m)}, \quad  n,m \geq 1.
\end{align*}

Let us now consider special cases of $p$ (below, $\delta$ denotes the dirac $\delta$-distribution).

\medskip

\subsubsection{ Right-skewed $\alpha$-stable distributions: $\quad p(d\beta) = \delta(\beta - \alpha) d \beta, \quad \ 0 < \alpha < 1$ }

\begin{figure}[h]
\begin{center}
\hspace*{-.8in}
\includegraphics[scale=.65]{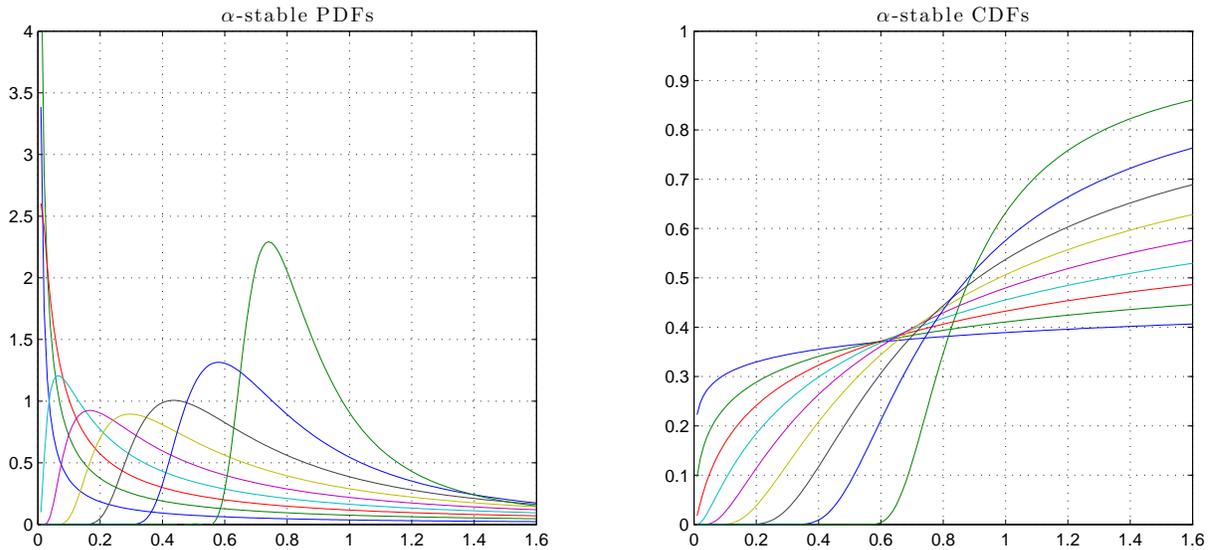}
\caption{ $\alpha$-stable PDFs (left) and CDFs (right) for $\alpha =
  i/10$ for $i=1,2,\dots,9$.  The PDF for  $\alpha = 0.9$ corresponds
  to  the right-most peak, and $\alpha = .1$ corresponds to the
  left-most peak (which isn't visible). The CDF with the steepest
  slope around $x = 0.8$ corresponds to $\alpha = 0.9$.    } \label{f:asPDFCDF}
\end{center}
\end{figure}

Indeed, an important example of
this distribution is the {\it right skewed $\alpha$-stable distributions}, for which
$p$ is a point mass at $\beta = \alpha$, with $0 < \alpha < 1$ (using  Proposition 1.2.11 in \cite{taqqu:1994}, this distribution
corresponds to the stable distribution $S(\alpha,1,\cos(\pi
\alpha/2)^{1/\alpha},0)$).  
These distributions lie in the family of scaling limits for sums of non-negative i.i.d.~random variables with infinite mean.  Form (\ref{e:asphin}), $\phi^{(n)}$ can be computed in various ways:
\begin{equation*}
\phi^{(n)}(\lambda) = (-1)^{n+1} \frac{\lambda^{\alpha - n} \Gamma(n -
  \alpha)}{\Gamma(-\alpha)} = - \lambda^{\alpha - n}
\prod_{m=0}^{n-1} (\alpha-m) = - \lambda^{\alpha - n} \sum_{m=1}^n \alpha^m S_n^{(m)}  
\end{equation*}
The PDFs and CDFs of the right skewed $\alpha$-stable distribution are
plotted using our method in Figure \ref{f:asPDFCDF} for several values
of $\alpha$.   For example, we compute the PDF of a $(1/2)$-stable
distribution at $x=1$ and obtain (in Mathematica)
\begin{equation*}
f(1) = 0.219695644733861
\end{equation*} 
This is exact to 15 decimal places and the computation took
approximately half a second.  An alternative method for this case is given in \cite{nolan:1997}.

\medskip

\subsubsection{Sums of right-skewed $\alpha$-stable distributions:
  $\displaystyle \quad p(d\beta) = \sum_{j=1}^r d_j \delta(\beta -
  \alpha_j) d \beta, \  d_j \geq 0, \  \displaystyle \ \sum_{j=1}^r d_j = 1$}

\begin{figure}[h]
\begin{center}
\hspace*{-.8in}
\includegraphics[scale=.65]{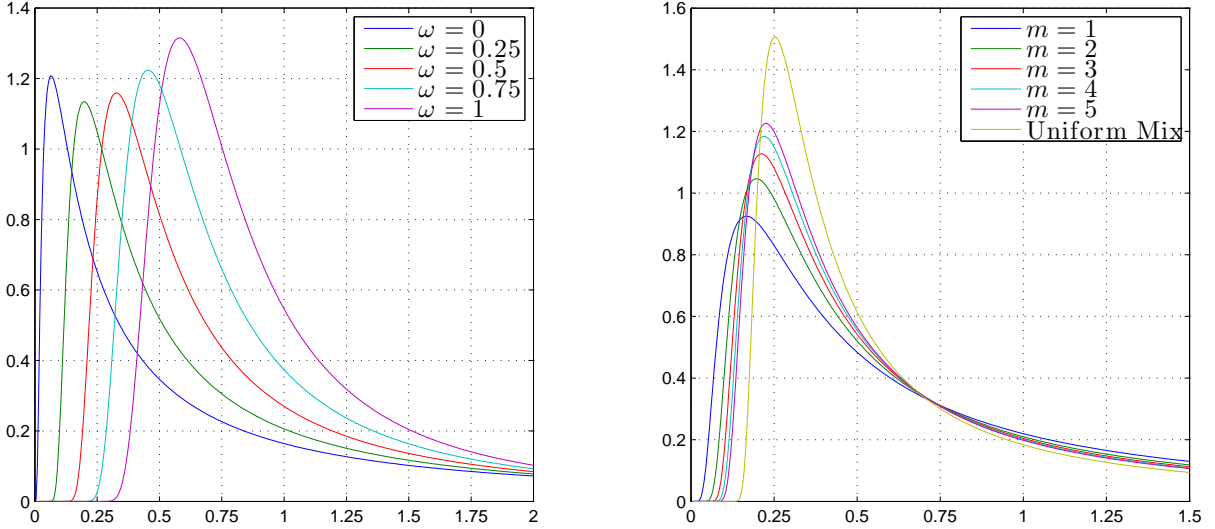}
\caption{ (Left figure) This is the plot of the PDF corresponding to the
  choice $p_\omega(du) = (1-\omega) \delta(x - 0.4) + \omega \delta(x - 0.8)$
for various choices of $0 \leq \omega \leq 1$. By increasing $\omega$,
the PDF for $\alpha = 0.4$ (on the left)  is morphing into the PDF for $\alpha =
0.8$ (on the right).  (Right figure)  These are PDFs corresponding to the sums $p_m(du) = \frac{1}{m}
\sum_{i=1}^{m} \delta(x - \frac{i}{m+1})$ for $m=1,2,\dots,5$ together
with the uniform continuous mixture on $(0,1)$, which has the highest peak.} \label{f:sumasPDFCDF}
\end{center}
\end{figure}

The previous case  can easily be generalized to weighted sums of independent $\alpha$-stable random variables.  From (\ref{e:asphin}) it follows that in this case,
\begin{equation*}  
\phi^{(n)} (\lambda) = (-1)^{n+1} \sum_{j = 1}^r d_j
\frac{\lambda^{\alpha_j - n} \Gamma(n - \alpha)}{\Gamma(-\alpha_j)} =-
\sum_{j=1}^r d_j \lambda^{\alpha_j - n} \prod_{m=0}^{n-1} (\alpha_j-m)
= - \sum_{j=1}^r d_j \lambda^{\alpha_j - n} \left( \sum_{m=1}^n \alpha_j^m S_n^{(m)}  \right).
\end{equation*}
Plot of PDFs of such distributions can be see in Figure \ref{f:sumasPDFCDF}.

\medskip

\subsubsection{ A ``uniform mixture'' of $\alpha$-stable
  distributions: $\quad p(d\beta) = d\beta, \ \beta \in (0,1)$ } \label{s:uniformmix}

This is an example of a distribution with no finite moments.  The
Laplace exponent (\ref{e:st-lap}) is given by 
\begin{equation*}
\phi(\lambda) = \int_0^1 \lambda^\beta d\beta =  \int_0^1 e^{\beta \log
  \lambda} d\beta =\frac{ (\lambda - 1)}{\log \lambda} \quad \mbox{if}
\ \ \lambda \neq 1,
\end{equation*}
and is $1$ if $\lambda = 1$. Higher derivatives of $\phi$ can be computed using  (\ref{e:asphin}).  The coefficients $c_m(\lambda)$, $m=1,2,\dots,n$ can be computed in a few different ways.  First, if $\log\lambda \leq 0$, we have
\begin{equation}\label{e:gaminc}
c_m(\lambda) = \int_0^1 \beta^m \lambda^\beta d \beta = \int_0^1
\beta^m e^{\beta \log \lambda} d \beta =\frac{ \gamma(1 + m,-\log\lambda)}{(-\log \lambda)^{m+1}}
\end{equation}
where $\gamma(a,b) $ is the lower incomplete gamma function (see above).  The formula (\ref{e:gaminc}) can still be used for $\log\lambda > 0$, however this requires analytic continuation of $\gamma$ which is not always easy to compute.

As another approach, we compute $c_m(\lambda)$ by treating the cases $|\log\lambda|<1$ and $\geq 1$ separately.  If $|\log\lambda|<1$ notice that
\begin{equation}
c_m(\lambda) = \int_0^1 \beta^m e^{\log(\lambda) \beta} d\beta = \sum_{j=0}^\infty  (\log\lambda)^j \int_0^1 \frac{\beta^{m+j} }{j!} d\beta = \sum_{j=0}^\infty \frac{(\log \lambda)^j}{(m+j+1) j! }.
\end{equation}
If $|\log\lambda| < 1$, taking the first 18 terms in this series gives
an absolute error of $< 10^{-15}$.   For $|\log\lambda| \geq 1$,  the $c_m$'s can be computed recursively by applying integration by parts:
\begin{eqnarray*}
c_0(\lambda) &=& \int_0^1 e^{\log(\lambda) \beta} d\beta =  \frac{\lambda - 1}{\log(\lambda)} \\
c_m(\lambda) &=& \int_0^1 \beta^m e^{\log(\lambda) \beta} d\beta = \displaystyle \frac{\lambda-  m c_{m-1}(\lambda)}{\log(\lambda)}, \quad m=1,\dots,n.
\end{eqnarray*}
Since this procedure involves a division by $\log \lambda$, it should
only be used when $|\log\lambda| \geq 1$.

\subsection{Integrals of non-random functions with respect to a Poisson random measure}

In this section we consider the large class of non-negative ID
distributions which can be expressed in terms of a Poisson stochastic
integral of a non random kernel.   For a recent review of these
integrals, see \cite{taqqu:2010}.

Let $(\Omega,{\cal F})$ be a probability space let ${\cal B}$ denote the Borel sigma field on $\mathbb{R}^n$.  Let $N(\cdot)$ denote a independently scattered Poisson random measure with control measure $\mu$, that is, $\mu$ is a measure on ${\cal B}$ and $N$ is a function $\Omega \times {\cal B} \rightarrow \mathbb{Z}^+ \cup \{0\}$ such that 
\begin{itemize}

\item[(i)] $N(A) \perp N(B)$ if $A,B \in {\cal B}$ are disjoint.
\item[(ii)] $\displaystyle N(\bigcup_{m=1}^n A_i) = \sum_{m=1}^n N(A_i)$ if $A_i$ are disjoint.
\item[(iii)] For each $A \in {\cal B}$, $N(A)$ has a Poisson distribution with mean $\mu(A)$.

\end{itemize}

Given such a pair $(N,\mu)$ one can define the following stochastic integral
\begin{equation} \label{e:poisint}
I(g) = \int_{\mathbb{R}^n} g( \xv ) N(d\xv)
\end{equation}
for a function $g$ on $\mathbb{R}^n$, for which $\int_{\mathbb{R}^n} \min(1,g(\xv)) \mu(d\xv) < \infty$ and, for our purposes, is non-negative.  In this case, the random variable $I(g)$ is also non-negative and has Laplace transform
\begin{equation}\label{e:PoisLT1}
\mathbb{E}e^{-\lambda I(g)} = \exp\left( - \int_{\mathbb{R}^n} (1 - e^{-\lambda g(\mathbf{s}) }) \mu(d\mathbf{s})  \right), \quad \lambda > 0.
\end{equation}

In order to compute the PDF and CDF of $I(g)$ using our method, (\ref{e:PoisLT1}) must first be rewritten in LK form.  In many cases, this can be done with some suitable change of variables  $\mathbf{u} = \Phi(\mathbf{s})$ satisfying $u_1 = (\Phi(\mathbf{s}))_1 =  g(\mathbf{s})$.  In this case, (\ref{e:PoisLT1}) becomes
\begin{equation}
\mathbb{E}e^{-\lambda I(g)}  = \exp\left(- \int_0^\infty (1 - e^{-\lambda u_1}) \Pi_g(d u_1) \right)
\end{equation}
where for any Borel set $A \subset \mathbb{R}^+$, the L\'{e}vy measure $\Pi_g(A)$ can be expressed formally as 
\begin{equation}\label{e:chofvar}
\Pi_g(A) = \int_{A \cap \mathbb{R}^{n-1}} |J(u_1,u_2,\dots,u_n)|   (\mu \circ \Phi^{-1} )(du_1,du_2,\dots,du_n),  
\end{equation}
where $|J|$ is the Jacobian $\partial \Phi^{-1}(u_1,\dots,u_n)/\partial (u_1,\dots,u_n)$. 

To illustrate this, let us now focus on special cases, in dimensions $n=1$ and $n=2$, where this change of variables can be made and our method applied.    

\medskip

\subsubsection{One dimensional Poisson integral}\label{s:oupoisson}

\medskip

\begin{figure}[h]
\begin{center}
\hspace*{-.8in}
\includegraphics[scale=.65]{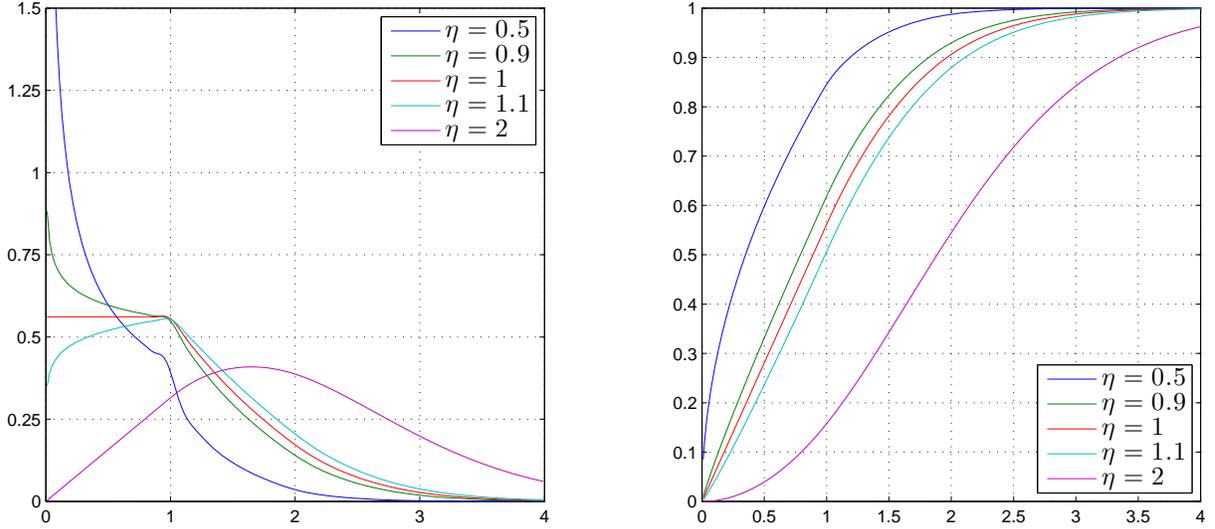}
\caption{ Plots of the PDF and CDF for the Poisson stochastic integral
  $I(e^{-s/\eta})$  with various values of the parameter
  $\eta$, which plays the role of a shape parameter.  The PDF for
  $\eta = 0.5$ is the highest on the left.    \label{f:OUpoisPDFCDF}  } 
\end{center}
\end{figure}

Assume $n=1$ and that the integrand $g$ is a monotone, non-negative function with inverse $g^{-1}$.  In this case, (\ref{e:PoisLT1}) can be rewritten in LK form using the change of variables $u = g(s)$:
\begin{equation*}
\mathbb{E}e^{-\lambda I(g)} = \exp\left( - \int_0^\infty (1-e^{-\lambda u}) \Pi_g(du) \right).
\end{equation*}
Since $u = g^{-1}(g(u))$, we get  $1 = (g' \circ g^{-1}) (u) \
(g^{-1})'(u) $ and hence,
\begin{equation}\label{e:newPi1}
\Pi_g(du)= \frac{\bone_{g((0,\infty))}(u)}{ |(g' \circ g^{-1})(u)|} (\mu
\circ g^{-1})(du) ),
\end{equation}
where $g((0,\infty))$ denotes the image of $(0,\infty)$ under $g$.
Suppose that we want to get the PDF and CDF of 
\begin{equation*}
I(g) = \int_0^\infty e^{-s/\eta} N(ds)
\end{equation*}
where $g(s) = e^{-s/\eta}$, $\eta >0$ is a parameter and the control measure $\mu$ is Lebesgue.  In this case, a simple calculation shows that (\ref{e:newPi1}) becomes
\begin{equation*}
\Pi_g(du) =  \left( \frac{\eta}{u} \right)  \bone_{(0,1]}(u)  du.
\end{equation*}
Thus, for this example,
\begin{equation*}
\phi(\lambda) = \int_0^1 (1-e^{-\lambda u})  \left(\frac{\eta}{u}\right) du = \eta \mathrm{Ein}(\lambda),
\end{equation*}
where $\mathrm{Ein}(\lambda)$ is the entire exponential integral defined
earlier. 
The derivatives $\phi^{(n)}(\lambda)$, for $n \geq 1$ can be given in closed form
\begin{equation*}
\phi^{(n)}(\lambda) = (-1)^{n+1} \eta \int_0^1 u^{n-1} e^{-\lambda u}
du = \frac{ (-1)^{n+1} \eta }{\lambda^n} \gamma(n,\lambda) = \frac{
  (-1)^{n+1} \eta }{\lambda^n} \left( (n-1)! - e^{-\lambda}
  \sum_{m=0}^{n-1} \frac{(n-1)!}{m!} \lambda^m \right),
\end{equation*}
where the last equality follows from (\ref{e:lincgam}).  We've plotted
the PDF and CDF of the random variable $I(e^{-s/\eta})$ for various
values of $\eta$ in Figure (\ref{f:OUpoisPDFCDF}) using our method.
Note that since the range of integration here is $(0,\infty)$, the
method described in \cite{veillette:2009a} doesn't readily apply.

\subsubsection{Integration with respect to non-negative L\'{e}vy
  process with a non-negative kernel}\label{s:ougamma}

In this example, we generalize the previous case by looking at
integration with respect to a non-negative L\'{e}vy process, or
equivalently, a one-dimensional non-negative ID random measure $L$
with control measure $\mu$.  $L$ is a random measure which satisfies the same conditions as the Poisson random measure $N$, except condition $(\mathrm{iii})$ is replaced by
\begin{itemize}
\item[(iii)'] There exists a L\'{e}vy measure $\Pi$ such that for any $A \in {\cal B}$, the distribution of $L(A)$ has Laplace transform
\begin{equation*}
\mathbb{E} e^{-\lambda L(A)} = \exp \left( - \mu(A) \int_0^\infty (1 - e^{-\lambda u}) \Pi(du)  \right)   .
\end{equation*}
\end{itemize}
Notice that the Poisson random measure $N$ corresponds to the choice $\Pi(du)  = \delta(u - 1)du$.  For a given function $g$, the stochastic integral $I_L(g)$ can be defined in terms of a two-dimensional Poisson stochastic integral:
\begin{equation}\label{e:ougamma}
I_L(g) = \int_0^\infty g(s) L(ds) \equiv \int_0^\infty \int_0^\infty u g(z) N(du,dz),
\end{equation}
where the control measure of $N$ is given now by $\Pi(du) \mu(dz)$.  Observe that the kernel $g$ must now satisfy $\int_0^\infty \int_0^\infty \min(1,u g(z))\Pi(du
) \mu(dz) < \infty$

Assume now for simplicity that $\Pi(du) = \Pi(u) du$ for some function
$\Pi$ and that $g$ is a non-negative monotone function with inverse $g^{-1}$.  In this case, we can obtain the LK form corresponding to $I_L$:
\begin{eqnarray*}
\mathbb{E}e^{-\lambda I_L(g)} &=& \exp \left( -\int_0^\infty
  \int_0^\infty (1 - e^{-\lambda u g(z) }) \Pi'(u) du \mu(dz) \right) \\
&=& \exp \left( - \int_0^\infty (1 - e^{-\lambda v }) \Pi'_g(v) dv \right)
\end{eqnarray*}
where we have made the change of variables $v = u g(z)$, $z' = z$ and the measure $\Pi_g(dv)$ is given by
\begin{equation}\label{e:newPi2}
\Pi'_g(v) = \int_0^\infty \Pi'\left( \frac{v}{g(z')} \right) \frac{1}{g(z')} \mu(dz').
\end{equation} 

To demonstrate our method in this case, consider 
\begin{equation*}
I_L(g) = \int_0^\infty e^{-z/\eta} L(dz) = \int_0^\infty \int_0^\infty
u e^{-z/\eta} N(du,dz)
\end{equation*}
where $g(z) =
e^{-z/\eta}$ with $\eta > 0$, $\mu$ Lebesgue, and $\Pi(du) = \kappa u^{-1}
e^{-u/\theta} du$, which is the L\'{e}vy measure corresponding to the
Gamma distribution with shape $\kappa>0$ and scale $\theta>0$
(\cite{Applebaum:2004}, Example 1.3.22).   With the change of variables $w = u e^{z/\eta}$, (\ref{e:newPi2}) implies
\begin{equation*}
\Pi'_g(u) = \int_0^\infty \Pi'\left( u e^{z/\eta} \right) e^{z/\eta} dz
= \frac{\eta \kappa}{u} \int_u^\infty w^{-1} e^{-w/\theta} dw =
\frac{\eta \kappa}{u} \Gamma\left( 0,\frac{u}{\theta} \right)
\end{equation*}
We can now compute the corresponding Laplace exponent for this case:
\begin{align}
\phi(\lambda) = \int_0^\infty (1 - e^{-\lambda u}) \Pi'_g(u) du &=
\eta \kappa \int_0^\infty \frac{(1 - e^{\lambda u} )}{u}
\Gamma\left(0,\frac{u}{\theta} \right) du \nonumber \\
&= \eta \kappa \int_0^\infty  \int_0^\lambda  e^{-u t}
\int_{u}^\infty 
\left( w^{-1} e^{- w/\theta} \right) dw \ dt \ du \nonumber\\
&=\eta \kappa \int_0^\lambda \left[ \int_0^\infty \int_0^{w}   ( w^{-1} e^{-u t -
    w/\theta} ) du \ dw \right] dt  \nonumber\\
&= \eta \kappa \int_0^\lambda \frac{1}{t} \left[ \int_0^\infty w^{-1}
  e^{-w/\theta} ( 1- e^{-w t} ) dw \right] dt \label{e:almostdilog}\\
&= \eta \kappa \int_0^\lambda \frac{1}{t}  \int_0^\infty \left(
  e^{-w/\theta} \int_0^t e^{-w s} ds \right) dw dt \nonumber\\
&= \eta \kappa \int_0^\lambda \frac{1}{t} \left[ \int_0^t
  \int_0^\infty e^{-w ( s + 1/\theta) } dw \ ds \right] dt \nonumber\\
& = \eta \kappa \int_0^\lambda \frac{ \log( 1 + t \theta) }{t} dt =
\eta \kappa L_2 ( 1 + \lambda \theta )  \label{e:dilog}
\end{align}
The derivatives of $\phi$ can also be computed exactly in this case.
Using (\ref{e:derofphi}) and (\ref{e:lincgam}),
\begin{align*}
\phi^{(n)}(\lambda) &= (-1)^{n+1} \eta \kappa  \int_0^\infty   u^{n-1}
e^{-\lambda u} \Gamma\left( 0 , \frac{u}{\theta} \right) du  \\
&= (-1)^{n+1} \eta \kappa \int_0^\infty u^{n-1} e^{-\lambda u} \int_u^\infty ( w^{-1} e^{-w/\theta} ) dw du \\
&=  (-1)^{n+1} \eta \kappa \int_0^\infty  w^{-1} e^{-w/\theta} \left[
  \int_0^w u^{n-1} e^{-\lambda u} du \right] dw \\
&=  (-1)^{n+1} \eta \kappa \int_0^\infty  w^{-1} e^{-w/\theta} \left[
  \frac{(n-1)!}{\lambda^n} - e^{-w \lambda} \sum_{m=0}^{n-1}
  \frac{(n-1)!}{m!} \frac{w^m}{\lambda^{n-m}} \right] dw \\
&=(-1)^{n+1} \eta \kappa  \frac{(n-1)!}{\lambda^n} \left[ \int_0^\infty w^{-1} e^{-w/\theta}
  ( 1- e^{-w/ \theta}) dw - \sum_{m=1}^{n-1} \frac{\lambda^m}{m!}
  \int_0^\infty w^{m-1} e^{-w(\lambda + 1/\theta)} dw \right] \\
&=(-1)^{n+1} \eta \kappa  \frac{(n-1)!}{\lambda^n} \left[ \log(1 + \lambda  \theta) -
  \sum_{m=1}^{n-1} \frac{ (\lambda \theta)^m }{m \ (1 + \lambda  \theta)^m} \right],
\end{align*}
where the first integral in the second to last line above is computed
as in (\ref{e:almostdilog}).   Alternatively, $\phi^{(n)}(\lambda)$ above can also be given simply in terms
of the Gauss hypergeometric function $\
_{2}F_1(\alpha,\beta,\gamma,x)$  (\cite{gradshteyn:2007}, Equation
6.455.1, page 657)
\begin{equation*}
\phi^{(n)}(\lambda) =(-1)^{n+1} \eta \kappa  \frac{\theta^n (n-1)!}{n} \ _{2}F_1(n,n,n+1,\lambda \theta). 
\end{equation*}

\begin{figure}[h]
\begin{center}
\hspace*{-.8in}
\includegraphics[scale=.65]{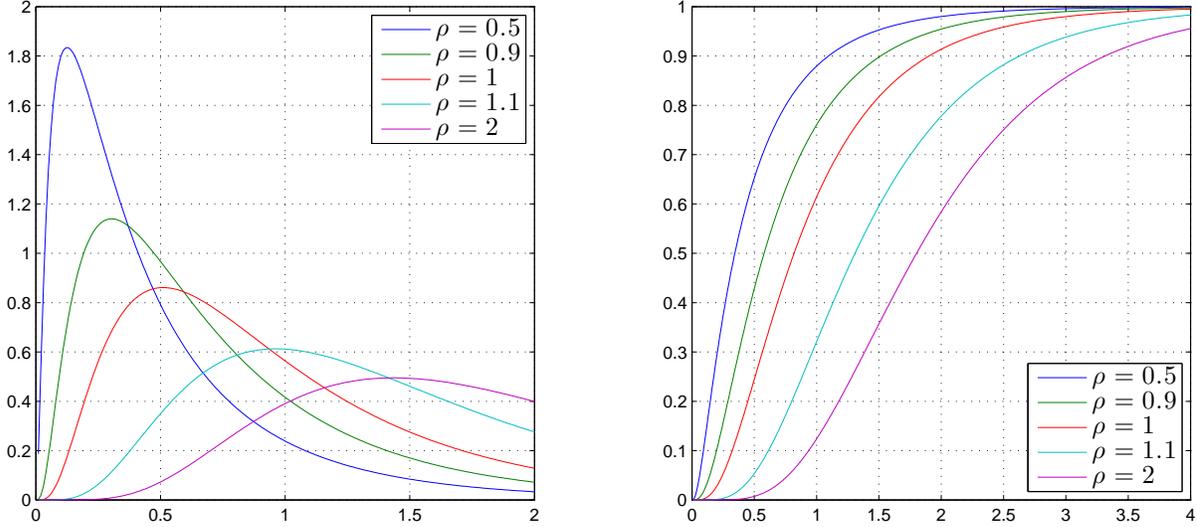}
\caption{ Plots of the PDF and CDF for the Poisson stochastic integral $I_L(e^{-s/\eta})$  with various values of the parameter $\eta$.  } \label{f:ougamma}
\end{center}
\end{figure}

In Figure \ref{f:ougamma} we've plotted the PDF and CDF of (\ref{e:ougamma}) for $\theta = 1$ and various values of the product $ \rho = \eta \kappa$ using our
method.

\section{Guide to software}\label{s:guidetosoft}

In this section, we will explain how to use the software written to
implement the method discussed in this paper.  Versions
of this code exist in MATLAB and Mathematica, and are freely available
by request form the authors.  Each version will include a file
containing examples to assist in using the code.

To begin using the code, download the file \texttt{NNINFDIV.zip} and extract the
directory \texttt{NNINFDIV}.  This directory contains both the MATLAB
and Mathematica programs in separate folders.  We will now focus on
these separately.

\subsection{MATLAB version}

To use the MATLAB version, launch MATLAB and add the directory
\texttt{NNINFDIV/MATLAB} to MATLAB's default path by typing

\medskip

\noindent $\gg$ \texttt{path(path,`mypath/NNINFDIV/MATLAB')}

\medskip
\noindent where `\texttt{mypath}' is the path which leads to the
directory \texttt{NNINFDIV}.   You are now ready to use the code
provided in this package.

The main function is called \texttt{nninfdiv}.  This function takes in
5 arguments in the following order:
\begin{align*}
& \texttt{X} &\qquad  &\mbox{Scalar or vector of input values.  Must be
  positive.} \\
& \texttt{DIST} &\qquad &\mbox{A cell array which specifies the distribution and
  parameters (see below)} \\
& \texttt{FUNC} &\qquad &\mbox{Type of function: `\texttt{pdf}' or `\texttt{cdf}'} \\
& \texttt{METHOD} &\qquad &\mbox{Extrapolation method: `\texttt{polynomial}' or
  `\texttt{rational}'} \\
& \texttt{TOL} & \qquad & \mbox{Target relative error tolerance}
\end{align*}
The last three inputs \texttt{FUNC,METHOD} and \texttt{TOL} are optional,
and take default values \texttt{'pdf', 'polynomial'} and $10^{-6}$.
 The input \texttt{DIST} is a cell array which contains the name of the
desired distribution followed by the parameters.  Possibilities for
\texttt{DIST} are
\begin{align*}
&\texttt{ \{'chi-squared',df,1\} } & \quad &\mbox{Chi-squared
  distribution with \texttt{df} degrees of freedom.} \\
&\texttt{ \{'chi-squared',df,[c1,...,cn]\} } & \quad &\mbox{Weighted
  sum of chi-squared distributions} \\
& \texttt{ \{'alpha stable',a,c\} } &\quad &\mbox{Alpha-stable
  distribution with $\alpha$ = a and scaling $c$}. \\
& \texttt{ \{'alpha stable',[a1,...,an],[c1,...,cn]\} } &\quad
&\mbox{Sum of $n$ weighted alpha-stable
  distributions}. \\
&\mbox{\texttt{ \{'uniform mix'\} }} & \quad & \mbox{Uniform mix from Section
  \ref{s:uniformmix}.} \\
&\mbox{\texttt{ \{'ou poisson',eta\}} } &\quad & \mbox{The integral
  $I(e^{-s/\eta})$ from Section \ref{s:oupoisson} } \\
&\mbox{\texttt{ \{'ou gamma',eta,kappa\}}} & \quad & \mbox{The integral
$I_L(e^{-s/\eta})$ from Section \ref{s:ougamma}}.
\end{align*}
The scaling constants seen in the alpha-stable and chi-squared
examples above compute the PDF/CDF of the scaled random variables $c
X$ for $c > 0$ in the single alpha-stable and single chi-squared case.
Likewise, in the weighted case, they return the PDF/CDF of $c_1
X_1 + c_2 X_2 + \dots c_n X_n$ with $c_i > 0$, and $X_i$ is
chi-squared with \texttt{df} degrees of freedom or is alpha-stable
with $\alpha =$\texttt{ai}.

Since the input for \texttt{nninfdiv} is long, it is often useful to
define a {\it function handle} in order to call the function more
easily.  For example, consider the $\alpha$-stable distribution with
$\alpha = 2/3$.  We define the PDF of this distribution in the
variable \texttt{f} by typing

\medskip

\noindent $\gg$ \texttt{f = @(x) nninfdiv(x,\{'alpha stable',2/3,1\},'pdf','polynomial',1e-6); }

\medskip

\noindent The function \texttt{f} now computes the PDF of the
alpha-stable distribution  with $\alpha = 2/3$ to within a relative error  of $10^{-6}$ using the
polynomial interpolation method.  For example, you may now type 
\begin{align*}
&\gg \texttt{f(1)}  &\qquad &\mbox{Computes the PDF at $x=1$}\\
&\gg \texttt{f([1 2 3])}  &\qquad &\mbox{Computes the PDF at $x=1,2$
and $3$.}\\
&\gg \texttt{plot([.05:.05:2],f([.05:.05:2]))} & \qquad &\mbox{Plots the PDF
  on the interval (0,2]}
\end{align*}

\medskip

\noindent {\bf Remark:} Obtaining relative errors less than $10^{-6}$
is sometimes difficult.  If your error tolerance cannot be reached, the program will return the
best estimate possible in double precision.  If high precision is preferred over speed, the Mathematica
version should be used.

\subsection{Mathematica version}

To use the Mathematica version, launch \texttt{Mathematica}
open the file \texttt{NNINFDIV.nb} located in the directory
\texttt{NNINFDIV/Mathematica}.  Once this file is open, select all its
contents by pressing \texttt{alt-a} on a PC or \texttt{cmd-a} on a
Mac.  Then compile the code by pressing \texttt{shift-return}.   You
are now ready to use the code in a separate notebook.

The main program is called \texttt{NNInfDiv} (capitalization
matters).  This program is called with 4 arguments:
\begin{align*}
& \texttt{X} &\qquad  &\mbox{Input value.  Must be a
  positive scalar.} \\
& \texttt{DIST} &\qquad &\mbox{A list which specifies the distribution and
  parameters (see below)} \\
& \texttt{FUNC} &\qquad &\mbox{Type of function: "\texttt{PDF}" or "\texttt{CDF}"} \\
& \texttt{TOL} & \qquad & \mbox{Relative error tolerance}
\end{align*}
The last two inputs \texttt{FUNC} and \texttt{TOL} are optional,
taking default values ``\texttt{PDF}'' and $10^{-15}$ respectively.
Possibilities for \texttt{DIST} include
\begin{align*}
&\texttt{ \{"Chi-Squared",\{c1,..cn\}\} } & \quad &\mbox{Sum of weighted
  chi-squared with weights \texttt{c1,...,cn}.} \\
& \texttt{ \{"Alpha Stable",a,c\} } &\quad &\mbox{Alpha-stable
  distribution with $\alpha$ = \texttt{a} and scaling \texttt{c}}. \\
&\mbox{\texttt{ \{"Uniform Mix"\} }} & \quad & \mbox{Uniform mix from Section
  \ref{s:uniformmix}.} \\
&\mbox{\texttt{ \{"OU Poisson",eta\}} } &\quad & \mbox{The integral
  $I(e^{-s/\eta})$ from Section \ref{s:oupoisson} } \\
&\mbox{\texttt{ \{"OU Gamma",eta,kappa\}}} & \quad & \mbox{The integral
$I_L(e^{-s/\eta})$ from Section \ref{s:ougamma}}.
\end{align*}
The scaling constants seen in the chi-squared and alpha-stable cases
above refer to the random variables $c X$ with $c>0$ in the
alpha-stable case and $c_1 X_1 + \dots c_n X_n$ in the chi-squared
case, with $X_i$ i.i.d chi-squared.

To simplify the call to this function, one can make a user defined
function.  For example, to make a function \texttt{F} which computes the CDF of
an $\alpha$-stable distribution with $\alpha = 2/3$, one can type

\medskip

\noindent \texttt{F[x\_] := NNInfDiv[x,\{ "Alpha Stable" , 2/3 , 1 \}
  , "CDF" ]  }

\medskip
\noindent The function \texttt{F} now computes the CDF of
$\alpha$-stable distribution with $\alpha = 2/3$ to a relative
precision of $10^{-15}$.  For example, one can now enter
\begin{align*}
&\texttt{F[1]}  &\qquad &\mbox{Computes the CDF at $x=1$}\\
&\texttt{Table[F[x],\{x,\{1,2,3\}\}]} &\qquad &\mbox{Computes the CDF
  at $x=1,2$ and $3$} \\
&\texttt{ListPlot[Table[\{x,F[x]\},\{x,0,2,.05\}],Joined -> True]} & \qquad &\mbox{Plots the CDF
  on the interval [0,2]}
\end{align*}

\medskip

\noindent {\bf Remark:}  Using  \texttt{NNInfDiv} with Mathematica's
\texttt{Plot} function is very slow, which is why we used the
\texttt{ListPlot} function above.  For faster plotting and function evaluation, the MATLAB
version of the code should be used.

\noindent Mark Veillette (mveillet@bu.edu) \& \\
\noindent Murad Taqqu (murad@math.bu.edu) \\
\noindent \small Dept. of Mathematics \\
\noindent \small Boston University \\
\noindent \small 111 Cummington St. \\
\noindent \small Boston, MA 02215


\end{document}